\input amstex
\documentstyle{amsppt}
\input psfig
\NoBlackBoxes
\magnification=1200
\hcorrection{.25in}
\advance\vsize-.75in

\define\biv{\Bbb B^4}

\define\hiii{\Bbb H^3}
\define\hiv{\Bbb H^4}
\define\hn{\Bbb H^n}
\define\riii{\Bbb R^3}
\define\isom{\operatorname{Isom}}
\define\volm{\operatorname{Vol}}
\define\intr{\operatorname{int}}

\define\const{4\pi^2/3}
\define\fai{\Phi_1}
\define\faii{\Phi_2}

\topmatter
\title
Finite-volume hyperbolic 4-manifolds that share a fundamental polyhedron
\endtitle

\author Dubravko Ivan\v si\' c
\endauthor

\address 
University of Oklahoma,
Mathematics Department,
Norman, OK 73019-0315
\endaddress

\email divansic\@aardvark.ucs.ou.edu
\endemail

\rightheadtext{Finite-volume hyperbolic 4-manifolds \dots}

\thanks
This paper represents a part of the author's PhD thesis.  The author would
like to thank his advisor, Dr. Boris Apanasov, for many helpful
conversations and suggestions concerning these results.
Research at MSRI is supported in part by NSF grant DMS-9022140.
\endthanks


\abstract 
It is known that the volume function for hyperbolic manifolds of
dimension $\geq 3$ is finite-to-one.  We show that the number of 
nonhomeomorphic hyperbolic 4-manifolds with the same volume can be made
arbitrarily large.  This is done by constructing a sequence of finite-sided
finite-volume polyhedra with side-pairings that yield manifolds.
In fact, we show
that arbitrarily many nonhomeomorphic hyperbolic 4-manifolds may
share a fundamental
polyhedron.  As a by-product of our examples, we also show in a
constructive way that the set of volumes of hyperbolic 4-manifolds
contains the set of even integral multiples of $\const$.   This is
"half" the set of possible values for volumes, which is the integral
multiples of $\const$ due to the Gauss-Bonnet formula
$\volm (M)=\const\cdot\chi(M)$.
\endabstract
\endtopmatter

\document

\head
0. Introduction and statement of results
\endhead

The original aim of research that produced this paper was to construct
non-compact hyperbolic 4-manifolds by means of side-pairings of
polyhedra.  Previous examples of hyperbolic manifolds with dimension higher
than
three were restricted to constructions via arithmetic
groups (see, for example, \cite{A2}, \cite{M}), or via "interbreeding" 
arithmetic groups to get non-arithmetic ones (\cite{G-P}) 
and there was only one (compact)
example using side-pairings, that of Davis in \cite{D}.
We were able to produce a number of examples of side-pairings of
hyperbolic 4-polyhedra and get, in addition, 
new information about volumes of hyperbolic 4-manifolds.
Further research led to consideration of embedability of
these manifolds as complements of surfaces in compact 4-manifolds ---
we deal with this in \cite{I}.

It is known (see \cite{W}) that for every constant $c>0$ there
are only finitely many complete non-homeomorphic hyperbolic $n$-manifolds
with volume $<c$, where $n\ge 4$.  For $n=3$ the set of volumes is a 
well-ordered (infinite) set, but still only finitely many manifolds may
have the same volume.  We concern ourselves with whether there is a bound
on the number of manifolds that have the same volume.

In dimension 3, this has been answered by Wielenberg (see \cite{Wi1})
for the non-compact case, and by Apanasov and Gutsul (\cite{A-G}) for
the compact one.  In both papers, for $N$'s that can be made arbitrarily
large, polyhedra are constructed in $\hiii$ and different side-pairings
are given on them whose quotient spaces are $N$ non-homeomorphic
hyperbolic manifolds.

In this paper we prove the analogous result for the non-compact case
in dimension 4, namely, 

\proclaim{Theorem A}
Given any number $N$, there exist more than $N$ non-homeomorphic,
non-compact, complete hyperbolic 4-manifolds of finite volume that share the
same fundamental polyhedron in $\hiv$.  In particular, they have the same
volume.
\endproclaim

The proof is by constructing polyhedra in $\hiv$ with different
side-pairings and utilizing Poincare's polyhedron theorem to see that
identifying paired sides yields complete hyperbolic $4$-manifolds. The
manifolds are then distinguished by how many ends they have.
It is known (Theorem 5.39 in \cite{A2} or \cite{A1}) that a
complete, hyperbolic,  geometrically finite $n$-manifold
has finitely many ends. If the manifold has finite volume, then
all of the ends are standard cusp ends,
that is, they are of the form $E\times [0,\infty)$, where $E$
is a closed flat manifold.  Furthermore, each end of the manifold corresponds
to a cycle (equivalence class under identification by side-pairing) of ideal
vertices of the polyhedron.
We count classes of ideal vertices for each of the side-pairings
that we construct and show that we arrive at different numbers for different
side-pairings.  Therefore, the resulting manifolds are non-homeomorphic,
because they have different numbers of ends.

We also give a geometric interpretation of the manifolds we
construct.  It turns out that each of the manifolds may be obtained
by taking two basic manifolds, cutting them along a two-sided totally geodesic
embedded 3-manifold and stringing several of these together by gluing
them along the isometric cuts.  In the process of justifying 
this interpretation
we prove a convenient sufficient condition for when a plane intersecting
a fundamental polyhedron for a group $G$ is precisely invariant with respect
to some subgroup $J\subset G$ (Theorem 4.4).

After the constructions in this work have been completed
the author became aware of two other preprints where
non-compact hyperbolic 4-manifolds were obtained 
by side-pairings of polyhedra in ways different from the one here.
Those constructions were then used to prove interesting results.
One of the preprints is \cite{N}, where B. Nimershiem constructs 
classes of examples that are used to show that the set of
all flat three manifolds that appear as cusps of hyperbolic
four-manifolds is dense in the set of all flat three manifolds.

The other, by J. Ratcliffe and S. Tschantz (\cite{R-T}), classifies all
non-compact hyperbolic 4-manifolds of minimal volume.
In addition to that, it is proved that the set of volumes is the positive
integral multiples of $\const$.  A by-product of our construction is

\proclaim{Theorem B} The set of all volumes of hyperbolic 4-manifolds
contains the even multiples of $\const$.
\endproclaim

This result is only "half as good" as the quoted one, but an advantage
is that we provide an explicit side-pairing to produce manifolds with the
desired volumes, whereas Ratcliffe and Tschantz's proof was not constructive
and gave only their existence.

\head
1. The polyhedron $P$ and its side-pairings $\fai$ and $\faii$
\endhead

We use the upper half-space model of hyperbolic space
to define a convex four-dimensional hyperbolic polyhedron
$P$ as an intersection of some hyperbolic half-spaces.

Recall that the Poincar\'e upper half-space model of hyperbolic $n$-space is
$\hn=\allowbreak\{(x_1,\dots,x_{n-1},t)\in{\Bbb R^n}\mid t>0\}$ with the metric
given by $ds^2={dx_1^2+\dots+dx_{n-1}^2+dt^2\over t^2}$.  The {\it
boundary at infinity} of a set $S$ is the set of all points in 
$\partial\hn={\Bbb R^{n-1}}\cup\{\infty\}$ that are in the
(Euclidean) closure of S. In the upper-half-space model
hyperbolic hyperplanes are either Euclidean half-spheres or Euclidean
half-planes orthogonal to $\partial\hn$
and they are uniquely determined by their own
boundaries at infinity, which are Euclidean $(n-2)$-spheres or $(n-2)$-planes
in ${\Bbb R^{n-1}}\cup\{\infty\}$.  
We will say that the hyperplane is {\it based},
respectively, on a sphere or a plane. (In our case $n=4$, so the hyperplanes
will be based on 2-spheres and 2-planes in $\riii$.) The angle
between hyperplanes is the same as the angle between their boundaries
at infinity.

Every hyperplane in $\hn$ determines two closed half-spaces: each
contains the hyperplane and their interiors are disjoint.  The 
$(n-2)$-sphere or -plane on which the hyperplane is based divides
$\Bbb R^{n-1}$ into two closed sets, each of which is the boundary of 
one of the half-spaces that the hyperplane determines. 
  
By a
{\it polyhedron}  in $\hn$ we will mean a connected subset of $\hn$ with
non-empty interior whose boundary is a locally finite collection of hyperplanes.
(The polyhedra in our construction are going to be intersections of
finitely many half-spaces, so they will also be convex.)
A {\it codimension-one side S} of $P$
is a subset of $\partial P$ such that $S=P\cap X$ and 
$S=\text{cl}_X(\text{int}_X S)$,
where $X$ is a hyperplane that bounds one of the defining
half-spaces of $P$.  Then $S$ is an $(n-1)$-dimensional convex polyhedron
in $X$.  Proceeding inductively we may define a {\it codimension-$i$ side}
of $P$ to be a codimension-one side of a codimension-$(i-1)$ side of
$P$.  (For more details on polyhedra, consult \cite{A2} or \cite{E-P}.)

Since every codimension-$i$ side is a polyhedron in
dimension $n-i$, we also call it an $(n-i)$-side.  Codimension-one sides
we will simply call {\it sides}, codimension-two sides we call {\it edges}, 
and we will use the term {\it vertex} for a $0$-side of $P$. Vertices of
$P$ are also called {\it finite vertices} or {\it real vertices} as opposed
to {\it vertices at infinity} or {\it ideal vertices} that are the
isolated boundary points of $P$ in $\partial\hn$.
To simplify notation, a hyperbolic
hyperplane, the side of $P$ lying on the hyperplane and the boundary
at infinity of the hyperplane will be denoted by the same letter.
(No confusion should arise here because our $P$'s are convex.)

\midinsert
$$\hss\psfig{file=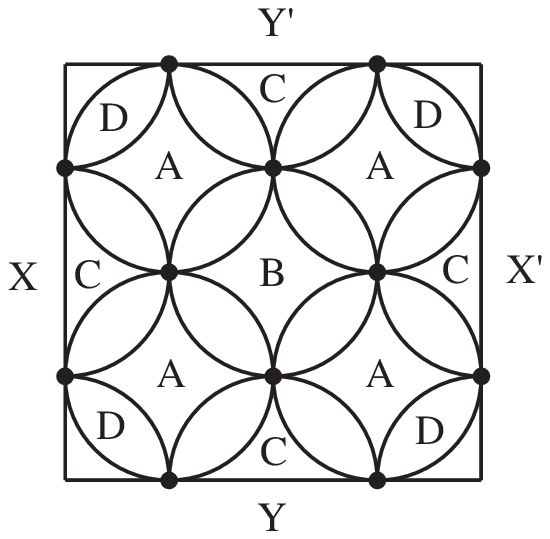,height=2truein}\hss$$
\botcaption{Figure 1}
Section of $P$ for $t=\sqrt2/2$, $z=\sqrt2/2$ showing
the real vertices
\endcaption
\endinsert

Consider the planes that bound the rectangular box $R\subset \riii$, 
$R=[-2,2]\times[-2,2]\times[-2\sqrt2,2\sqrt2]$.  Add to them the 12 spheres
of radius $\sqrt2$ with centers
$(\pm1,\pm1,j2\sqrt2)$ for $j=-1,0,1$ and the 18
spheres of the same radius with centers $(j,k,\pm\sqrt2)$ for $j,k=-2,0,2$.
The upper part of figure 2, going from left to right, depicts intersections
of these spheres with planes with constant $z$-coordinates 
$-2\sqrt2,-\sqrt2,0,\sqrt2,2\sqrt2$.
Label the spheres by the letters $A_i,A'_i,B_i,B'_i,C_i,C'_i,D_i,D'_i$
in either of the ways suggested by figure 2.
Let $X_1,X'_1,Y_1,Y'_1,Z_1,Z'_1$ be respectively the planes
$\{x=-2\},\ \{x=2\},\ \{y=-2\},\ \{y=2\},\ \{z=-2\sqrt2\},\ \{z=2\sqrt2\}$.

Each of the planes that comprise the boundary of $R$ and
each of the above spheres determine a hyperplane in 
$\hiv=\{(x,y,z,t)\in \Bbb R^4  \mid t>0\}$ that divides $\hiv$ into two
half-spaces.   
For the spheres we choose the half-spaces whose boundary at infinity is
unbounded in
$\riii$, for the planes the half-spaces so that the intersection of their
boundaries at infinity is the rectangular
box $R$. The polyhedron $P$ is defined as the intersection of those half-spaces.
For later convenience, we set $P_-=\{ (x,y,z,t)\in P\mid z\le0 \}$,
$P_+=\{ (x,y,z,t)\in P\mid z\ge0 \}$.

The following observations about the spheres and planes that we just defined
are easy to check.
\roster
\item Any two spheres that intersect do so at an angle of $\pi/2$.
\item Whenever the intersection is non-empty,
spheres $A_i,A_i'$ intersect planes $X_1,X'_1,Y_1,Y'_1$ at
angle $\pi/4$.  Any other pair of spheres or planes with non-empty
intersection intersects at angle $\pi/2$.
\item $R$ is completely covered by the closed balls bounded by the spheres.
This means that $P$ has finite volume and has only finitely many points in its
boundary at infinity.
\item $P$ has 36 vertices at infinity, which correspond to points not covered by
the open balls.  Their position is illustrated in Figure 9.
\endroster

It is not obvious right away that $P$ also has finite vertices (that is,
0-sides that are in $\hiv$).  This is because many sets of four hyperplanes
bounding the polyhedron $P$ meet at one point.  For example, sides
$A_1,B_1,C_1$ and $A_2'$
meet at the point $(0,1,-3\sqrt2/2,\sqrt2/2)\in\hiv$, and sides
$A_1,C_1,D_1$ and $Y_1'$ meet at the point $(-1,2,-3\sqrt2/2,\sqrt2/2)$.
Figure 1 depicts the section of $P$ where $t=\sqrt2/2$ and
$z=j\sqrt2/2$.  Here
$j$ is any of $-3,-1,1,3$  as the section for every $j$ is the same.  We can
see where four sides of $P$ intersect in a vertex and what letters those sides
are labeled by. Figure 1 shows the location of all the vertices in one
section --- there being $4\cdot 12=48$ in all.

\topinsert
$$\hss\psfig{file=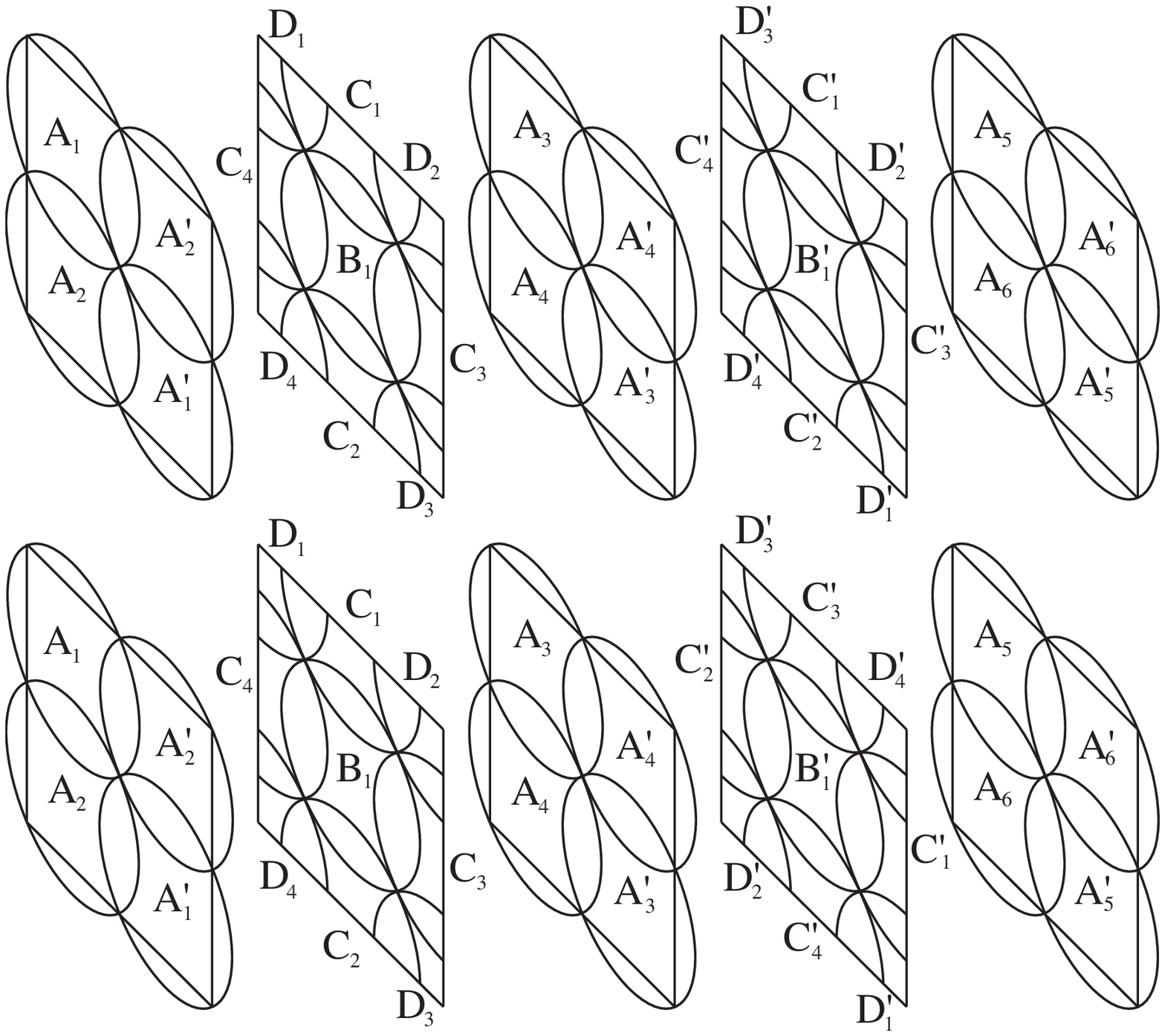,height=4truein}\hss$$
\botcaption{Figure 2}
Side-pairings $\fai$ (top) and $\faii$ (bottom)
\endcaption
\endinsert

Now we are ready to define two ways to pair sides of $P$. A 
{\it side-pairing $\Phi$} of $P$ is a rule which associates to each side
$S$ of $P$ a side $S'$ of $P$ and an isometry $s$ which sends $S$ to
$S'$.  This rule is subject to the conditions that
$s(\intr P)\cap\intr P=\emptyset$, that the side associated to $S'$
must be $S$, that the isometry that takes $S'$ to $S$ must be $s^{-1}$.
For more details, see \cite{E-P}
or \cite{R}.

First, define the
following isometries of $\riii$:  
$$\align
q_0&=\text{reflection in plane}\  \{z=0\} \\
q_1&=\text{reflection in plane}\  \{x-y=0\} \\
q_2&=\text{reflection in plane}\ \{x+y=0\} \\
s_1&=\text{rotation by $\pi$ about line}\ \{x+y=0, z=0\} \\
s_2&=\text{rotation by $\pi$ about line}\ \{x-y=0, z=0\} \\
t_0&=\text{translation by $2\sqrt2$ in the $z$ direction}.
\endalign
$$
Use the same letters to denote the extensions of these maps to $\hiv$.
(A Euclidean isometry $f:\riii\to\riii$ extends to a hyperbolic isometry
given by $(x,y,z,t)\mapsto(f(x,y,z),t)$.)
 
Let $i_S$ denote the reflection in the hyperplane $S$. By $s$ we denote the
hyperbolic isometry that pairs the sides $S$ and $S'$ (it sends $S$ to $S'$). 
We define $\fai$ to be the side-pairing given in the upper-half-space model
by
$$
\align
x_1 &=\text{translation by 4 in the $x$ direction} \\
y_1 &=\text{translation by 4 in the $y$ direction} \\
z_1 &=t_0^2=\text{translation by $4\sqrt2$ in the $z$ direction}\\
b_1 &=q_1\circ t_0 \circ i_{B_1}  \\
a_j &=q_l \circ i_{A_j},\ \text{so that}\ l\equiv j\ (mod\ 2) \\
c_k &=q_0 \circ i_{C_k} \\
d_k &=q_1 \circ t_0\circ i_{D_k},
\endalign
$$
where $j=1,\dots,6,\ k=1,\dots,4$ and $l=1,2$.  The upper half of Figure 2
shows which sides are paired.

To get another side-pairing, $\faii$, we alter $\fai$ in
the way the sides labeled by $B$'s, $C$'s and $D$'s are paired.  Refer to
the lower half of Figure 2 to see
which sides are paired.  We define the new pairings $b_1, c_k$ and $d_k$ by
$$
\align
b_1 &=q_0 \circ i_{B_1}  \\
c_k &=q_1 \circ t_0 \circ i_{C_k} \\
d_k &=s_l \circ q_l \circ i_{D_k},\ \text{so that}\ l\equiv k\ (mod\ 2),
\endalign
$$
where $k=1,\dots,4$ and $l=1,2$.

\head
2. Two closed finite-volume hyperbolic 4-manifolds 
\endhead

In this section we prove 

\proclaim{Theorem 2.1}
\roster
\item"{(i)}" The side-pairings $\fai$ and $\faii$ 
generate discrete torsion-free subgroups $G_1$ and $G_2$ of $\isom\hiv$ whose
fundamental polyhedron is $P$. Therefore, the quotient of $\hiv$ by the action
of either of the groups is a complete hyperbolic 4-manifold.
\item"{(ii)}" $\hiv/G_1$ has seven ends while $\hiv/G_2$ has eight.  In
particular, the two manifolds are not homeomorphic.
\endroster
\endproclaim

\demo{Proof} 
To prove assertion $(i)$, we use Poincar\'e's polyhedron theorem.  For
details, the reader can consult \cite{E-P} and \cite{R} which
were our main references, while we will shortly state the version that we are
going to use.  Other versions of the polyhedron theorem may be found in
\cite{A2} and \cite{Ma}.

\midinsert
$$\hss\psfig{file=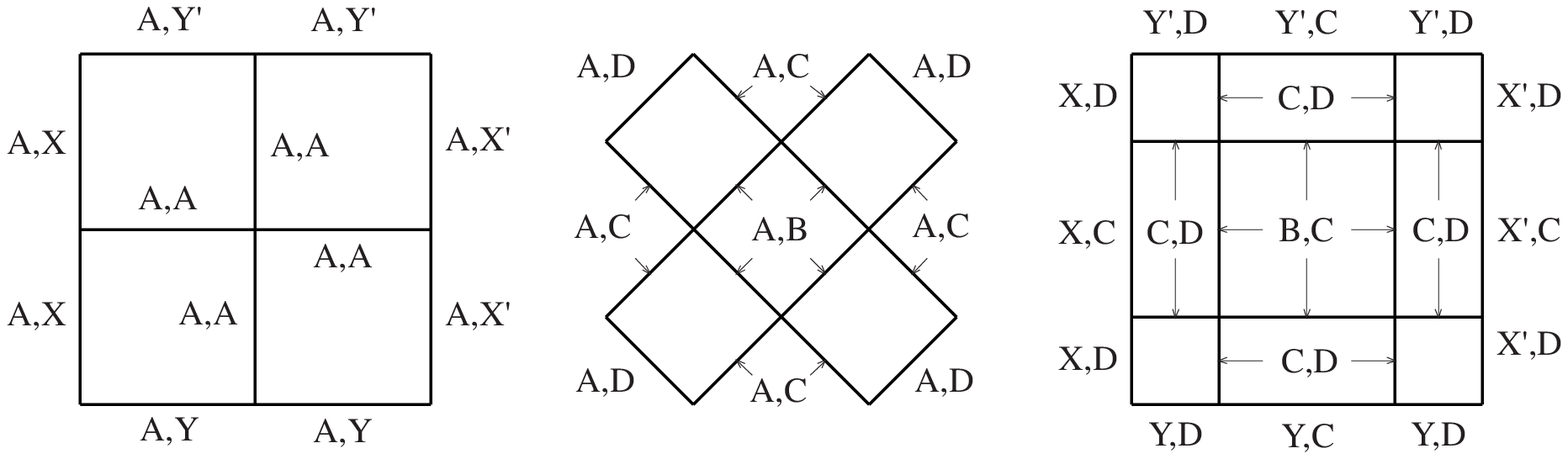,height=2truein}\hss$$
\botcaption{Figure 3}
Representing edges in $P$ by diagrams
\endcaption
\endinsert

First of all, the maps defined above really do map a side of $P$
isometrically  onto a side of $P$.  To verify this for $\fai$, 
notice that each of its side-pairings $a_1,\dots,a_6,c_1,\dots,c_4$,
is of the form 
$f\circ i_S$ where S is a side and $f$, {\it which preserves $P$}, is an
extension to $\hiv$ of a Euclidean transformation on $\riii$.
The $i_S$ keeps $S$ fixed so $f\circ i_S(S)=f(S)$, and
$f$, being an isometry of $P$, sends its sides to some other sides.  
We proceed similarly for the the side-pairings $b_1,d_1,\dots,d_4$:
each of them is of
the form $f\circ i_S$, but 
$f$ is now an isometry that takes $P_-$ to $P_+$, so it sends sides of
$P_-$ to sides of $P_+$.  But the sides $B_1,D_1,\dots,D_4$ are sides of
both $P_-$ and $P$, and, likewise the sides 
$B'_1,D'_1,\dots,D'_4$ are sides of both $P_+$ and $P$.
For the side-pairing $x_1$ ($y_1$ and $z_1$ are done similarly), compose it
with a reflection in the side $X_1'$ to get an isometry of $P$ whose image of
$X_1$ are the same as by $x_1$.  Therefore $x_1$
carries $X_1$ to another side of $P$, namely  $X_1'$.  
The claim is proved in the same way for the side-pairing $\faii$.

\remark{Remark 2.2} Recall that a horosphere in $\hn$ is either a
Euclidean sphere tangent to $\partial\hn$ or a Euclidean hyperplane
$\{t=c\}$ parallel to $\partial\hn$.  The former are said to be {\it
centered} at the point of tangency with $\partial\hn$, the latter at
$\infty$.
Consider a set $T$ of disjoint horospheres, each centered at a vertex
at infinity of $P$.  For the vertex $\infty$ choose, say,
the horizontal plane $\{t=3\}$. For the other vertices choose  
horospheres of the same radius that is small enough 
so that the horospheres intersect only those
sides of $P$ which contain the center of the horosphere.
If the center of a horosphere is on the boundary of a hyperplane $S$,
then $i_S$ preserves the horosphere.  This combined with an argument
like in the preceding paragraph can be used to show that the side
pairings $\fai$ and $\faii$ satisfy the 

\vskip6pt
{\it Consistent horosphere condition}: there exists a set $T$ of 
disjoint horospheres centered at 
ideal vertices of $P$ so that if $g$ is a side-pairing of
a side that contains the center of a horosphere $H\in T$ in its boundary, then
$g(H)$ is again a horosphere from $T$.

\endremark

\vskip6pt

Another condition for Poincar\' e's polyhedron theorem is the
"edge cycle condition", called Cyclic in 
\cite{E-P}.  In general, a side-pairing on $P$ induces an
equivalence relation on $P$ that is generated by the relation
$x\sim s(x)$, where $x\in\partial P\cap S$, $S$ is a side of $P$
and $s$ its side-pairing.  The equivalence class $[x]$ of $x$ under this
equivalence relation, is called the {\it cycle} of $x$.   The cycle of an
$i$-side is defined analogously, so that it contains all the $i$-sides
of $P$ that are identified by a string of side-pairings.

\midinsert
$$\hss\psfig{file=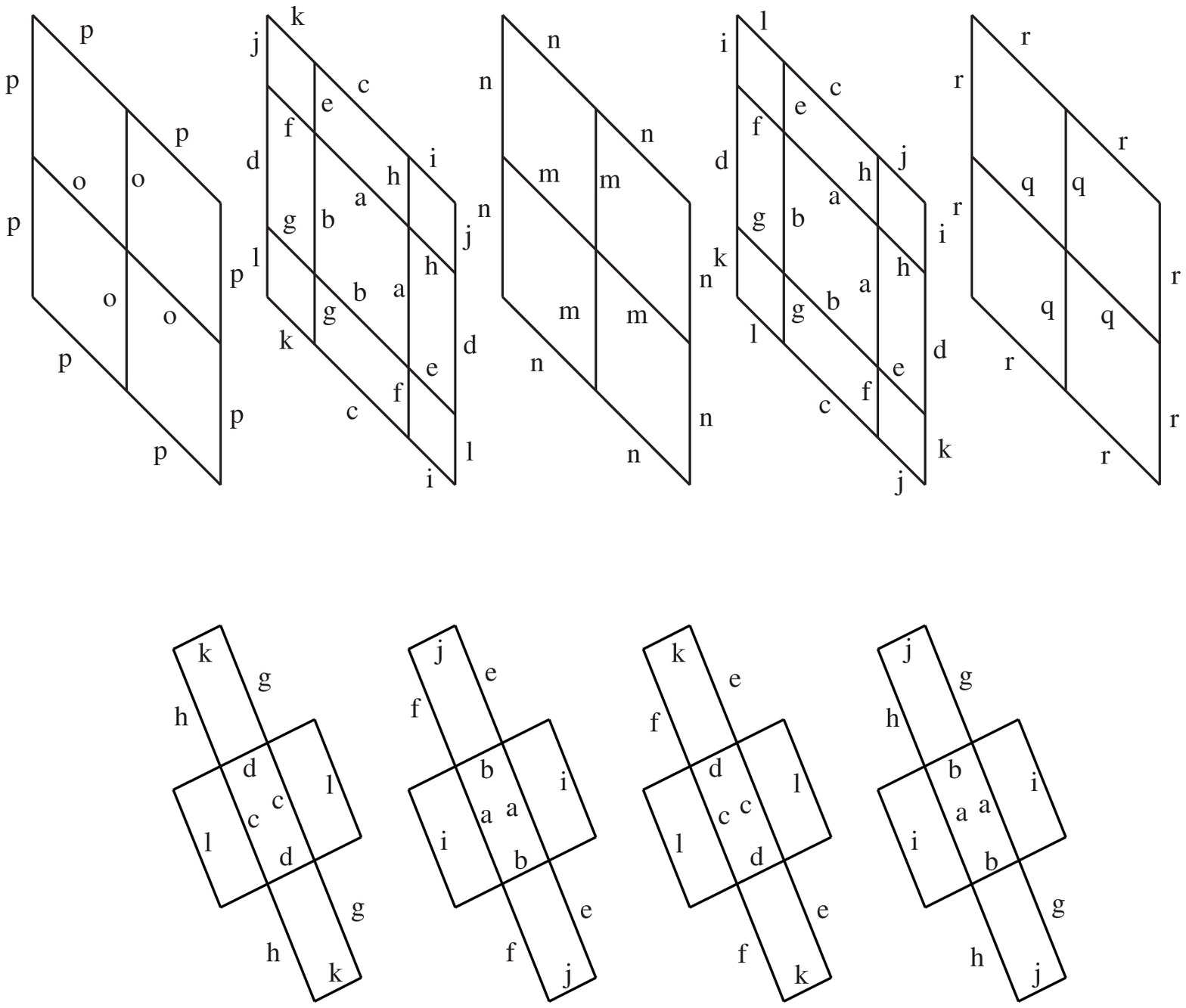,height=4truein}\hss$$
\botcaption{Figure 4}
Edge chase for $\fai$
\endcaption
\endinsert

Every edge (codimension-2 side) of $P$ is the intersection of two uniquely
determined
sides of $P$.  The {\it dihedral angle} at an edge is the angle in the
interior of $P$ that the two sides subtend.
A cycle of edges can be obtained in the following way.  Start with an
edge $E_1$, which is the intersection of sides $S_1$ and $R_1$, and let
$g_1$ be the isometry pairing $R_1$ and some side $S_2$ of $P$.  We get
that $g_1(E_1)=E_2$, where $E_2$ is an edge determined by $S_2$ and
some other side $R_2$.  Now let $g_2$ be the isometry pairing $R_2$
and some side $S_3$.  Continuing in the same way we get a sequence of edges,
sides and isometries $\{\sigma_i=(E_i,S_i,R_i,g_i)\}_{i=1,2\dots}$.
This procedure is commonly called "edge-chasing".
We require that the above sequence  have a period $q$ (called {\it first
cycle length in \cite{E-P}}), that is
$\sigma_{q+1}=\sigma_1$ for some $q$.  The cycle of edges will then
consist of exactly $E_1,\dots,E_q$.  Due to finite-sidedness of our $P$,
this condition will automatically be satisfied.

It is clear that $g_q\circ\dots\circ g_1(E_1)=E_1$,
but it may happen that the restriction of $g_q\circ\dots\circ g_1$
on $E_1$ is not the identity.  The second part of the edge cycle
condition is that there must be a $k$ so that
$(g_q\circ\dots\circ g_1)^k\mid _{E_1}=id$.  The number $kq$ is called
the {\it second cycle length} in \cite{E-P}.

Finally, to fulfill the edge cycle condition we must show that if
$\theta_i$, $i=1,\dots,q$ is the dihedral angle of edge $E_i$, then there
is a non-zero integer $m$ so that $k(\theta_1+\dots+\theta_q)=2\pi/m$.
We may now formulate Poincar\' e's polyhedron theorem for the case of a
hyperbolic polyhedron as follows.

\proclaim{Theorem 2.3} (Poincar\' e's polyhedron theorem)
Let $\Phi$ be a side-pairing on a polyhedron $P\subset\hn$ that satisfies
both the edge cycle condition and the consistent horosphere condition.
Then the side-parings of $\Phi$ generate a discrete group $G\subset\isom\hn$
whose fundamental polyhedron is $P$.
\qed\endproclaim

Now we check the edge cycle condition for edges of $P$ and the two
side-pairings $\fai$ and $\faii$.
With notation as above, we will always have $k=1$ and
$m=1$.  Hence, the second cycle length will always be the same as
the first cycle length
and they will be 4 and 8 respectively for edges with dihedral angles $\pi/2$
and $\pi/4$.  Therefore, the sum of dihedral angles will be exactly $2\pi$.

Firstly, we make sure that all edges are in cycles of said length. 
For all edges that are intersections of sides based on planes (i.e. "vertical"
sides), this
check is easy and boils down to checking the conditions of Poincar\' e's
polyhedron theorem for a rectangular Euclidean parallelepiped with parallel
sides paired by Euclidean
translations.  The check for any of the edges of type $A_i\cap Z_1$ is
also straightforward.  For all the other edges, we use the diagrams in
Figures 4 and 5 to simplify and visualize the task of verifying.

\midinsert
$$\hss\psfig{file=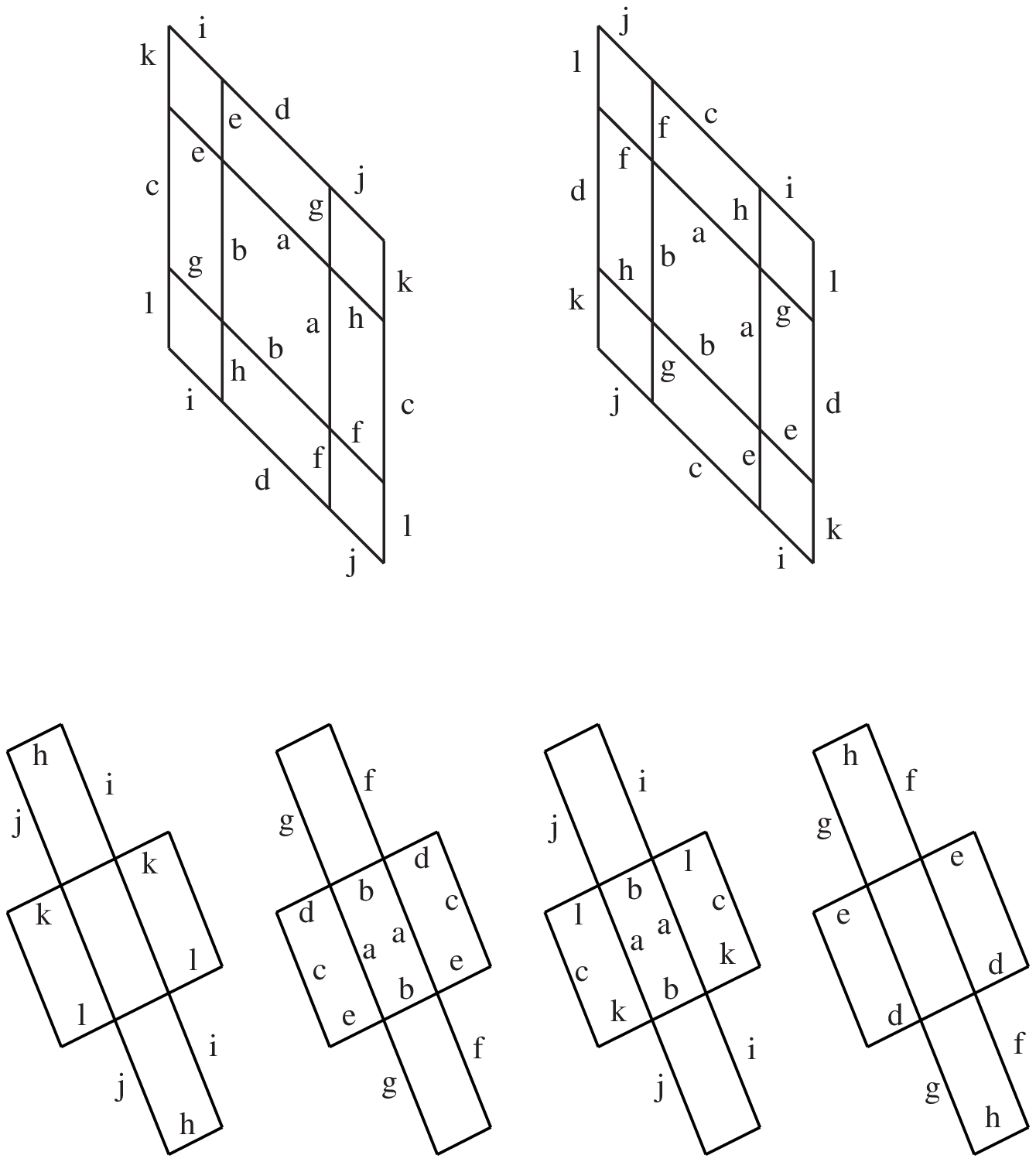,height=4truein}\hss$$
\botcaption{Figure 5}
Edge chase for $\faii$
\endcaption
\endinsert

Notice that the intersection of two hyperbolic hyperplanes is a codimension-2
hyperbolic subspace whose boundary at infinity is the intersection of
the boundaries at infinity of the hyperplanes.  Therefore, every edge
$E$ of $P$ lies on a codimension-2 subspace determined either by 
the intersection of a sphere and a plane or by the intersection of two
spheres in $\riii$.  (The spheres and planes are the boundaries at
infinity of the sides that determine $E$.)
These intersections are circles and they are 
represented by segments in Figure 3. The letters next to each segment indicate 
which side-types have generated the edge represented by it.
To get the left and right diagrams we take intersections of
planes and spheres labeled by the letters in the diagrams and then
project them to the plane $\{z=0\}$.  For the middle diagram, we first
project to the plane $\{z=0\}$ the centers of those pairs of spheres whose
labels are listed in it. Then we take the perpendicular bisector
of the line joining those centers. To account for all the edges,
we need several of these diagrams (Figure 4).		
Edge chasing for $\fai$ and $\faii$ is now performed on Figures 4 and 5 
respectively.  All edges in one cycle in each horizontal component
of the pictures are labeled with the same letter.  Some are not labeled
because their cycles are similar to other labeled cycles.  Also, Figure
5 omits some of the edges because their cycles are the same as for
$\fai$.

For example, choose the edge $A_1\cap A_2$.  The edge chase, yielding
the cycle labeled $o$ in upper part of Figure 4 is 

$$
A_1\cap A_2 @>a_2>> A_2'\cap A_1 @>a_1>> A_1'\cap A_2'
@>a_2^{-1}>> A_2\cap A_1' @>a_1^{-1}>> A_1\cap A_2
$$

As another example, choose the edge $A_1 \cap D_1$.  The edge chase, yielding
the cycle labeled $k$ in lower part of Figure 4 is
$$
A_1\cap D_1 @>d_1>> D_1'\cap A_3' @>a_3^{-1}>> A_3\cap D_3' 
@>d_3^{-1}>>D _3\cap A_1' @>a_1^{-1} >> A_1\cap D_1.
$$

Next, we check that for transformations $g_1,\dots,g_q$ 
obtained by edge-chasing we have 		
$g_q\circ\dots\circ g_1\vert_{E_1}$=1.  As before, $g_i=f_i\circ r_i$,
where $f_i$ is the extension to $\hiv$ of a Euclidean transformation on $\riii$
and $r_i$ is either reflection in a hyperplane containing $E_i$ or the 
identity. Let $f=f_q\circ\dots\circ f_1$. It is not difficult to see
that $f$ is always orientation preserving.
Clearly $g_q\circ\dots\circ g_1\vert_{E_1}=f\vert_{E_1}$, and it will be
enough to show that $f=1$.  We will need the following easy lemma.

\midinsert
$$\hss\psfig{file=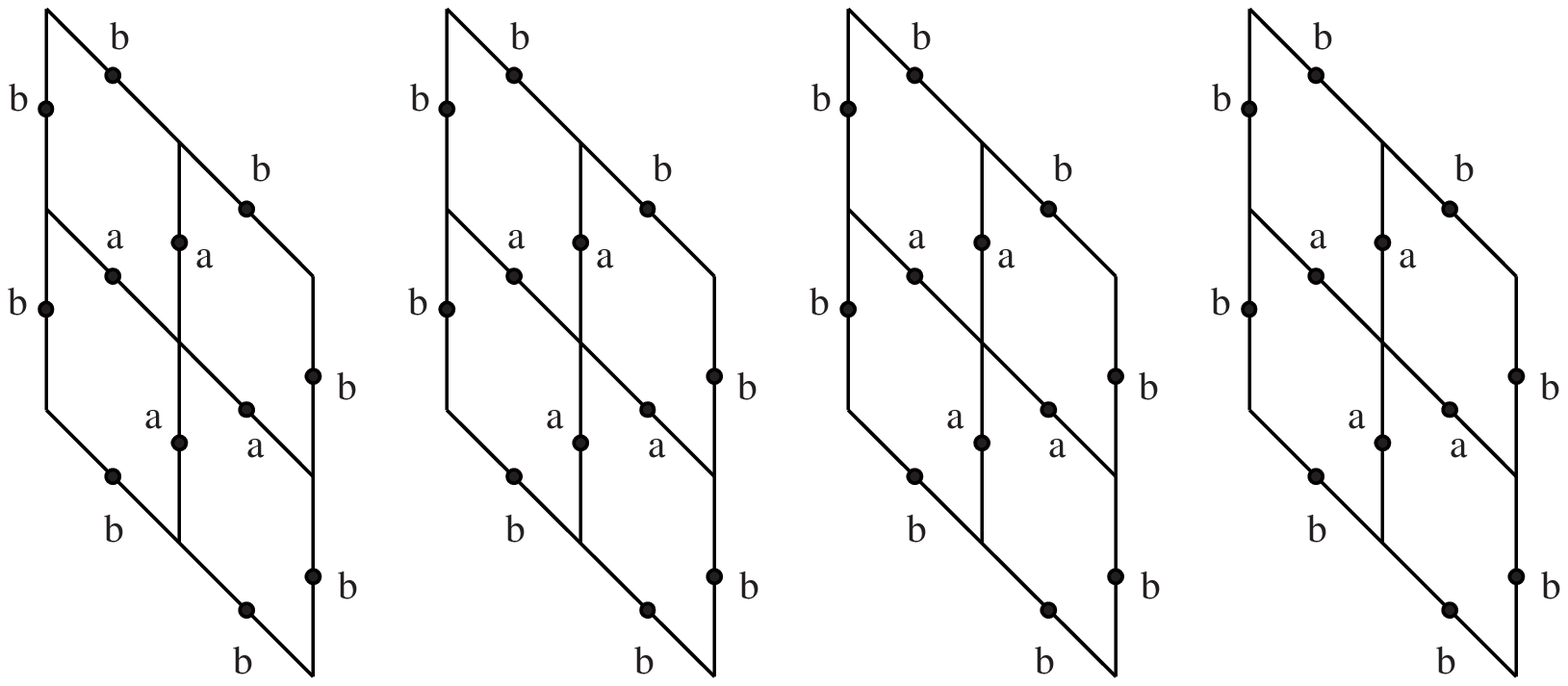,height=2truein}\hss$$
\botcaption{Figure 6}
Cycles of real vertices
\endcaption
\endinsert

\proclaim{Lemma 2.4}
Let $f$ be a nontrivial orientation-preserving Euclidean isometry of $\riii$
that preserves a circle.  Then it is a rotation about a line called the
{\rm axis} of $f$.
Moreover, we have:
\roster
\item"{(i)}" If we write $f$ as $Ux+u$, where $U$ is an orthogonal
transformation
and $u\in\riii$ then the axis of $f$ is parallel to the axis of $U$.  
\item"{(ii)}"The axis of $f$ passes through the center of the circle and is
either in the plane of the circle, or perpendicular to it.  In the first case,
the rotation is by angle $\pi$.
\item"{(iii)}" If $f$ preserves a line $l$, then its axis either orthogonally
intersects $l$ or it is that line.  
\qed\endroster
\endproclaim

Let $f=f_q\circ\dots\circ f_1$ as above and suppose it is nontrivial.
In what follows we interchangeably view $f$ as a Euclidean isometry on
$\riii$ or as a hyperbolic isometry of $\hiv$.

We know that $f$ preserves the circle that is the base of the edge $E=E_1$.
Also, f preserves the family of planes 
$V=\{(x,y,z)\in\riii\mid x=4k+2\ \text{or}\ y=4k+2,\ k\in \Bbb Z\}$ because
it is a composite of maps from $\{x_1, y_1, z_1, q_0, q_1, q_2, s_1, s_2\}$.
The rotational part $U$ of $f$ is
a composite of maps from $\{q_0, q_1, q_2, s_1, s_2\}$,
each of which preserves the axes $l_1$ and $l_2$ of
$s_1$ and $s_2$, so $U$ preserves them too.
Now, looking at possible
positions of the circle we get one of the following cases: 

{\it Case 1.} When $E$ is one of the sides represented in the lower half of
Figure 4 then, by $(i)$ and $(ii)$ of Lemma 2.4, the axis $l$ of $U$ must
either have direction vector $(\pm 1,\pm 1,\sqrt2)$ or is in the plane
perpendicular to that vector.  Since $U$ preserves $l_1$ and $l_2$, by part 
$(iii)$ of the
lemma, both of $l_1$ and $l_2$ must be either perpendicular or identical
to $l$.  If $l\perp l_1$ and $l\perp l_2$, then $l$ is the $z$-axis, which
contradicts the possible positions of $l$.  If $l$ is equal to either
$l_1$ or $l_2$, then one can see that the axis of $f$ is going pass exactly
through the segment that represents $E$ in the middle diagram of Figure 3.
However,
it is clear that no rotation about these segments can preserve the family
of planes $V$ so we must have $f=1$.

{\it Case 2.} When $E$ is one of the sides represented in the upper half of
Figure 4 the axis $l$ of $U$ lies, by parts $(i)$ and $(ii)$ of Lemma 2.4,
in one of the the planes $\{x=0\}$ or $\{y=0\}$ or it 
is the $x$- or $y$-axis.  Clearly $l\ne l_1$ and $l\ne l_2$, so applying
part $(iii)$ of Lemma 2.4 again we get $l\perp l_1$
and $l\perp l_2$, which means $l$ is the $z$-axis.   As long as $E$ is not
of the forms $C_i\cap X_1,\ C_i\cap Y_1,\ D_i\cap X_1$ or $D_i\cap Y_1$
the center of the circle on which $E$ is based has odd $x$ and $y$
coordinates.  However, a rotation about an axis parallel to $z$ through
such points cannot preserve the family of planes $V$ and we again get
$f=1$.

{\it Case 3.} In the remaining cases, if $E$ is of form $D_i\cap X_1$ or
$D_i\cap Y_1$, regard $f$ as an isometry of $\hiv$ and examine its
action on vertices of $P$ that are on $E$.  Looking at Figure 1, it is
clear that $f$, being by the above a rotation about an axis parallel to $z$,
must send vertices that are on $E$ to points whose either $x$ or $y$
coordinate falls out of $[-2,2]$, a contradiction with the fact that $f$
preserves $E$.  Finally, if $E$ is of form 
$C_i\cap X_1$ or $C_i\cap Y_1$ we just compute $f$: it is always $q_0^2=1$
for the side-pairing $\fai$ and it is always $q_1^2=1$ for the
side-pairing $\faii$.

Thus, we have shown that $f=1$ in all possible cases and the
edge cycle condition has been verified.

Since the consistent horosphere condition is fulfilled by Remark 2.2, we
may apply theorem 2.3 to get that the groups $G_1$ and $G_2$
generated by the side-pairings $\fai$ and $\faii$ are discrete, and that $P$ is
the fundamental polyhedron for both of them.

What we do not yet know is
whether $G_1$ and $G_2$ are torsion-free, that is, whether $\hiv/G_i,\ i=1,2$
are hyperbolic manifolds and not just orbifolds.

\midinsert
$$\hss\psfig{file=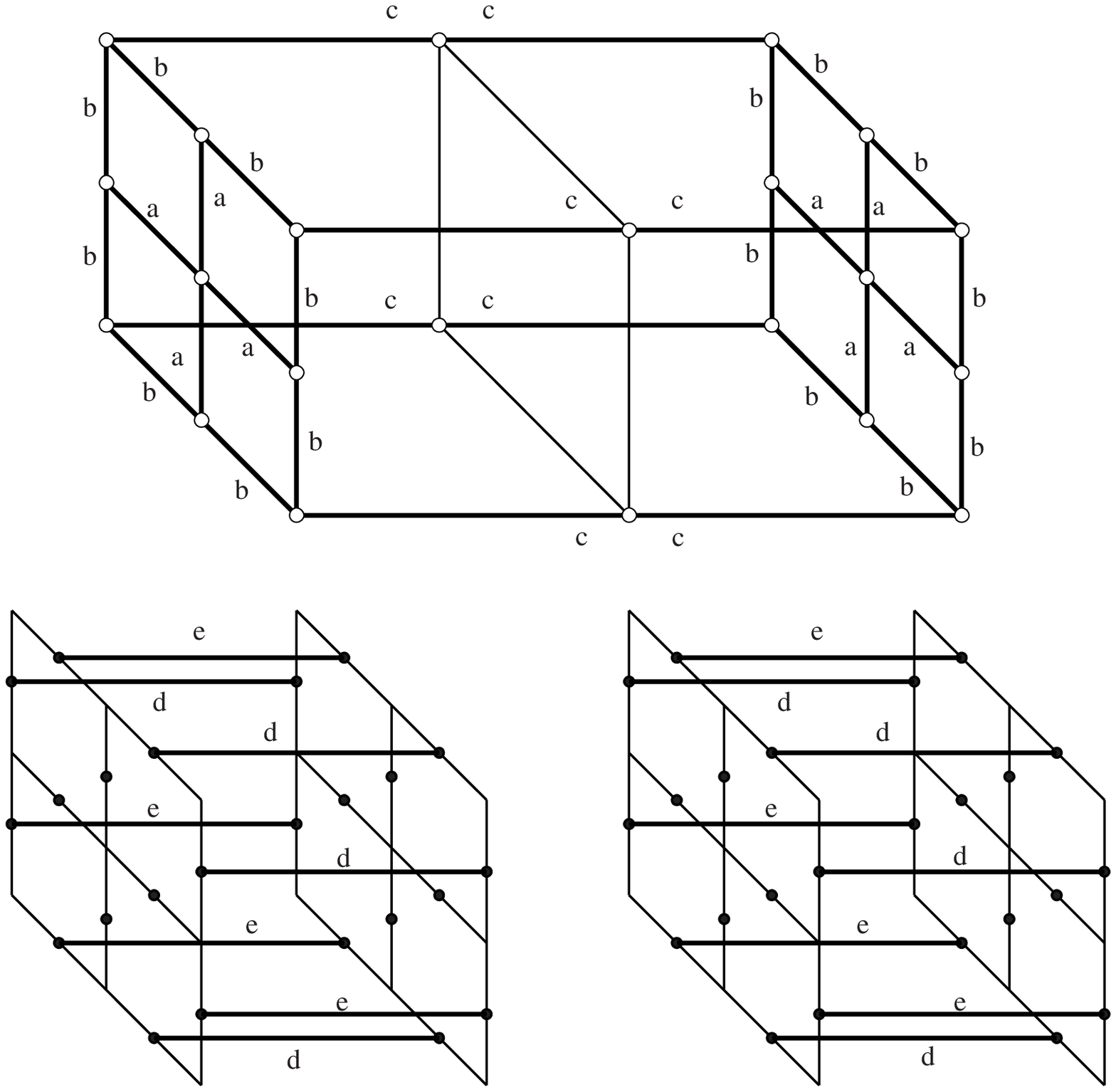,height=4truein}\hss$$
\botcaption{Figure 7}
Cycles of some 1-sides for $\fai$
\endcaption
\endinsert

Recall (see \cite{R}) that the {\it normalized solid angle at point $x$} of a
polyhedron $P$ is defined as $\omega(x)=(\volm B(x,r)\cap P)/ \volm B(x,r)$.
Here $B(x,r)$ is a hyperbolic ball about $x$ or radius $r$, $\volm$ is
hyperbolic volume and $r$ is taken
small enough so that  $B(x,r)$ intersects only
those sides of $P$ on which $x$ lies.  Let $[x]=\{x_1,\dots,x_n\}$
be the cycle of $x$ for some side-pairing of $P$. We define the
{\it normalized solid angle sum of
$[x]$} as $\omega[x]=\sum_{y\in[x]} \omega(y)$.  We are going to use
Theorem 11.1.1 from \cite{R} which says that

\proclaim{Theorem 2.5}
If $\omega[x]=1$ for every $x\in P$ then the group $G$ generated by the
side-pairings of $P$ is torsion-free.
\qed\endproclaim

Knowing that $P$ is a fundamental polyhedron for a discrete group $G$
implies $\omega[x]\le 1$ for every point of $P$.  Really, for every
$x_i\in[x]$, choose an isometry $g_i\in G$ taking $x_i$ to $x$.  (In
general, there may be many ways to make the choices.)  We now have a injective
map from $\{x_1,\dots,x_n\}$ to \{translates of $P$ under $G$ containing
$x$\} given by $x_i\mapsto g_i(P)$.  Since 
$\{g_1(P)\cap B(x,r),\dots,g_n(P)\cap B(x,r)\}$ fill
out maybe only a portion of $B(x,r)$, we get $\omega[x]\le 1$.
(Note that the strict inequality
will occur if and only if  $x$ is the fixed point of an element in $G$.)

Therefore, it is enough to see that $\omega[x]\ge 1$.
For an $x\in P$ that is in the interior of 3 or 4-sides of
$P$, it is clear that $\omega[x]=1$.  For an $x$ in the interior of 2-sides,
this is the edge-cycle condition.  This leaves 0- and 1-sides to be checked.

\midinsert
$$\hss\psfig{file=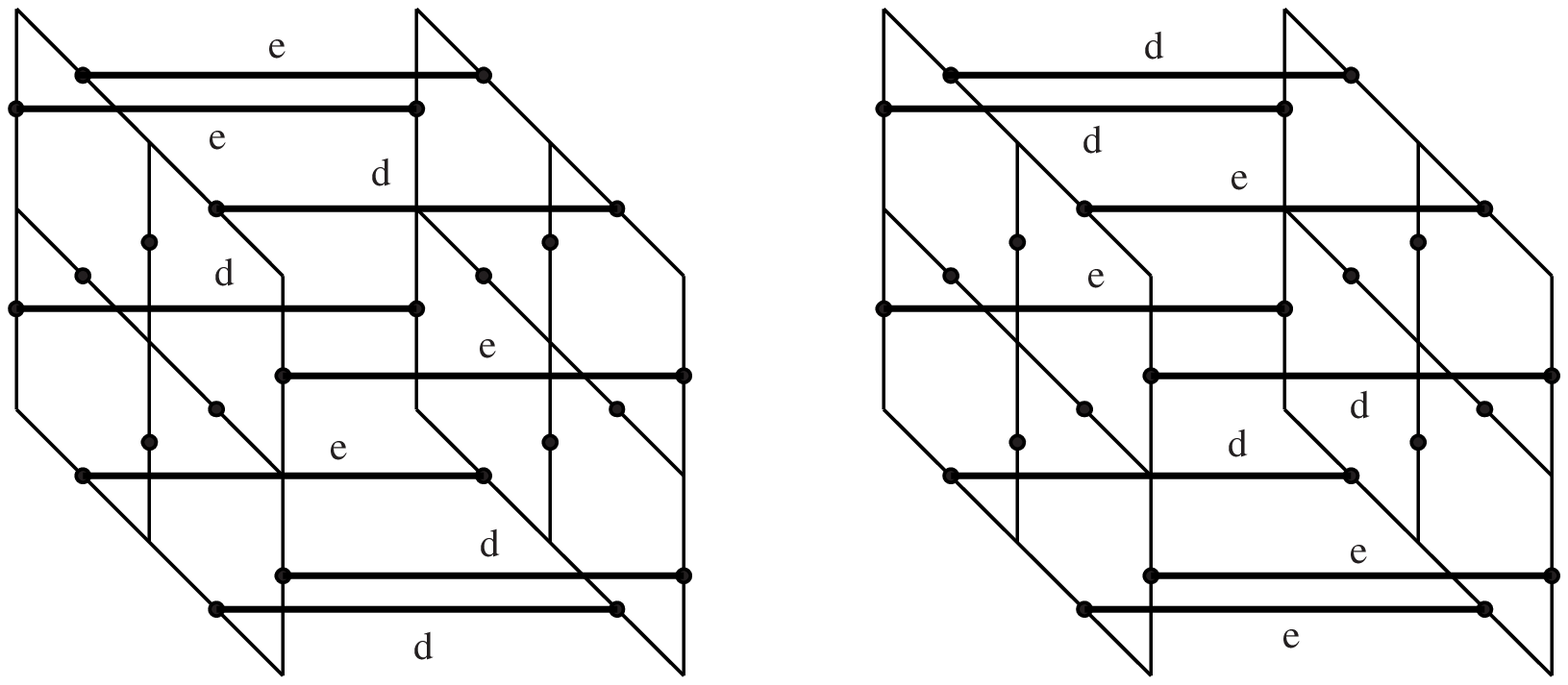,height=2truein}\hss$$
\botcaption{Figure 8}
Cycles of some 1-sides for $\faii$
\endcaption
\endinsert

Figure 6 shows the sections of $P$ for $t=\sqrt2/2$ and
$z=-3\sqrt2/2$,$-\sqrt2/2$,$\sqrt2/2$, $3\sqrt2/2$.  These sections contain
all the vertices (0-sides) of $P$.  All vertices in the
same cycle (there are only two cycles) are labeled by the same letter
--- the cycles are the same for
both $\fai$ and $\faii$. Vertices in the cycle labeled $a$ occur as the
intersection of four hyperplanes, each pair of which meets at angle $\pi/2$.
Normalize $P$ so that a vertex $x$ from cycle $a$ is the origin in the ball
model $\biv$ of hyperbolic space. We see that the normalized solid angle at $x$ 
is the same as the normalized solid angle at $0\in\biv$ subtended by the four
coordinate planes, and that is $1/16$.
Vertices in the cycle labeled $b$ are always intersections of four
hyperplanes where one pair of them intersects at angle $\pi/4$ and all other
pairs intersect at angle $\pi/2$.  Normalizing as before, we see that the
normalized solid angle at $x$ is the same as the one at $0$ subtended by
the hyperplanes $\{x_2=0\},\ \{x_3=0\},\ \{x_4=0\},\ \{x_1-x_2=0\}$, and that is
$1/32$.
Since cycles $a$ and $b$ contain 16 and 32 points respectively, we are done.

Now consider 1-sides.  There are three cases depending on whether a 1-side
$F$ connects an ideal and a real vertex, two ideal vertices, or two real
vertices.

For the first case, let $F$ be a 1-side of $P$ that is a geodesic
half-line between one real and one ideal vertex of $P$.  If for some
$x\in\operatorname{int} F$ we have $\omega[x]<1$, then $x$ is a fixed point of
some nontrivial $g\in G$. The isometry $g$ must preserve $F$ --- otherwise,
we'd have $g(F)\cap F=\{x\}$ and this contradicts the fact that translates
of $P$ meet only along $i$-sides. However, this implies that 
the real vertex on $F$ is fixed under $g$, a possibility we just proved
cannot happen.  So, we are left with checking 1-sides that have as endpoints
either both real or both ideal vertices of $P$.

Every 1-side is the intersection of three different sides of $P$. Hence,
to find all 1-sides with both endpoints real or ideal, we have to find
pairs of real or ideal vertices lying on the same three sides.  Figure 7	
schematically depicts those 1-sides of $P$.  A boldface line segment joining
the real or ideal vertices indicates the existence of a 1-side joining them.
The three
sides on which the 1-sides lie are easily deduced from their position in
the picture. (For example, the 1-sides labeled $c$ are intersections of sides
labeled by $D$'s, $X$'s and $Y$'s.)
As before, the letters on the 1-sides indicate to which
cycle of 1-sides they belong.

It takes a bit of checking to see that we have found all the needed
1-sides.  For example, to see that no 1-side joins an ideal vertex in the
plane $\{z=0\}$  to an ideal vertex in the plane $\{z=\sqrt 2\}$
 we note that every vertex in the plane $\{z=0\}$ lies on only one of
the $B,C$ or $D$-sides, and some $X,Y,Z$ or $A$-sides, while every
vertex in the plane $\{z=\sqrt2\}$ lies on only one of the $A$-sides and
some number of $B,C$ or $D$-sides.  Therefore, any pair of vertices
from those two hyperplanes cannot belong to the same three sides.

The 1-sides in the cycle labeled $a$ are intersections of three sides
meeting pairwise at angles $\pi/2$ (two sides labeled by  $A$ and one
by $Z$) .  Taking an $x$ from a side in the
cycle, and normalizing in $\biv$ so that $x=0$ and the three sides are
the first three coordinate hyperplanes, we see that $\omega(x)=1/8$.
Since there are 8 1-sides in the cycle, we get $\omega[x]\ge 1$.

Other cycles are checked in the same way.  For an $x$ on a 1-side in the
cycles $b,c,d,e$ we get normalized solid angles of  respectively 
$1/16,1/8,1/8,1/8$ with $16,8,8,8$ 1-sides in the cycle, so
$\omega[x]\ge 1$. Thus, we have proved that side-pairings $\fai$ and
$\faii$ give rise to hyperbolic 4-manifolds.

Assertion $(ii)$ of Theorem 2.1 is verified by counting cycles of ideal
vertices for $\fai$ and $\faii$.  The cycles are shown in Figure 9.
\qed\enddemo

\head
3. Proof of Theorem A
\endhead

\midinsert
$$\hss\psfig{file=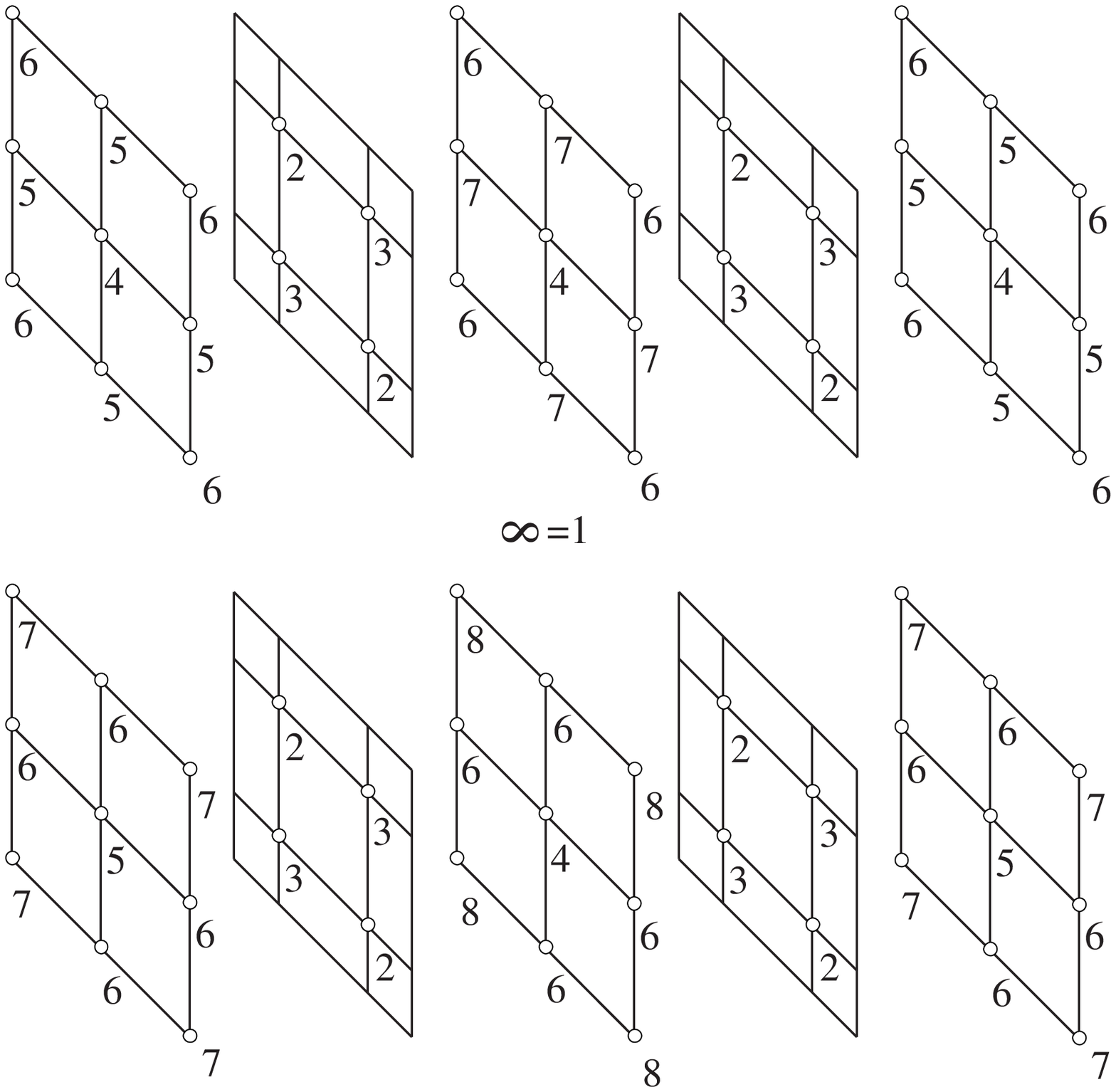,height=4truein}\hss$$
\botcaption{Figure 9}
Cycles of ideal vertices for $\fai$ (top) and $\faii$ (bottom). In both
cases, the vertex $\infty$ is labeled by 1
\endcaption
\endinsert

The crucial observation that the proof uses is that the side-pairings
$\fai$ and $\faii$ yield manifolds with respectively 7 and 8 ends.  This is
immediate from the fact that there are 7 and 8 classes of vertices at
infinity for the side-pairings $\fai$ and $\faii$, as shown in Figure 9.

Choose a positive integer $n$, and let $Q_1,\dots,Q_n$ be $n$ copies of
the polyhedron $P$.  Call each of them a {\it block}.  To each block,
assign either of the side-pairings $\fai$ or $\faii$ and call it a block 
{\it of type} $\fai$ or $\faii$.  Now form a new
polyhedron $Q$ by attaching side $Z_1$ of $Q_i$ to side $Z'_1$ 
of $Q_{i-1}$, $i=2,\dots,n$, i.e. by stringing the blocks together
in a linear fashion in the direction of the $z$-axis.  $Q$ has on it a
side-pairing induced by the side-pairings on each block.  Clearly, the sides 
that were attached have
vanished, so they don't fall under this rule: the remaining sides $Z_1$
of $Q_1$ and $Z'_1$ of $Q_n$ are paired by the translation
$t^{2n}$. A moment's reflection convinces us that the side-pairing on
$Q$ generates a torsion-free group $G$ --- the proof is basically an
$n$-fold repetition of the proof of the same result for the side-pairings
$\fai$ and $\faii$.
Except for $i$-sides, $i=0,1,2$, that are contained in the sides $Z_1'$ of
$Q_1$ and $Z_1$ of $Q_n$,
all the cycles are just inherited cycles from pairings on each block. 
 The special cases are
easily dealt with --- they follow patterns established for $\fai$ and
$\faii$.

How many cycles of vertices at infinity does $Q$ have?  Let there be $k$
blocks of type $\fai$ and $n-k$ blocks of type $\faii$ among
$Q_1,\dots,Q_n$.  Starting with $Q_1,\dots,Q_l$, adding a block of type 
$\fai$ onto $Q_l$ will
add three new cycles of vertices at infinity: there are seven cycles on
$\fai$ but four fall into cycles already existing on $Q_1,\dots,Q_l$.  
Similarly,
adding a block of type $\faii$ adds four new cycles of vertices at infinity.
It is now easy to see that $Q$ will have $4+3k+4(n-k)=4+4n-k$ cycles of
vertices at infinity.  We have complete freedom of choice for $k$, so by
varying $k$ from 0 to $n$ we can get manifolds with anywhere from $4+3n$
to $4+4n$ ends.  Therefore, we have obtained
at least $n+1$ nonhomeomorphic
manifolds with the same fundamental polyhedron $Q$, which is what we set out
to prove.  (We likely get many more, since we completely ignored the various
orderings of blocks of the two types that are possible when constructing
$Q$.)
\qed

\head
4. A geometric interpretation of the construction
\endhead

In this section we analyze the construction of the manifolds in the
previous section from a gluing-of-manifolds perspective.

Let $M=\hn/G$, where $G$ is a discrete torsion-free subgroup of
$\isom\hn$.  A {\it totally geodesic hypersurface} is a subset
$N\subset M$ so that for every $x,y\in N$ every geodesic connecting
$x$ and $y$ is also contained in $N$.  We are interested in embedded
totally geodesic hypersurfaces which are the ones for which
$p^{-1}(N)$ is a disjoint union of hyperplanes in
$\hn$, where $p:\hn\to M$ is the standard projection. Let $H$ be one of
those hyperplanes, and $J\subset G$ its stabilizer in $G$, that is the 
subgroup $J=G_H=\{j\in G\mid j(H)=H\}$. 
Then $H$ is {\it precisely invariant under} $J$, i.e. $g(H)=H$ when
$g\in J$ and $g(H)\cap H=\emptyset$ when $g\in G\setminus J$.
In particular, we want to look at some hypersurfaces with 
$\volm N<\infty$, so they will correspond to subgroups $G_H\subset G$
that act on a hyperplane $H$ as a hyperbolic lattice.

Conversely, we may start with a subgroup $J\subset G$ and a hyperplane $H$
precisely invariant under $J$.  Then $H/J$ is an embedded totally geodesic
hypersurface in $M$.

Let $M_i=\hiv/G_i$, 
where $G_i$ is generated by the side-pairing $\Phi_i,\ i=1,2$, and let
$H$ be the supporting hyperplane of side $Z_1$ of the polyhedron $P$.
Suppose we know that $H$ is precisely invariant under the subgroup 
$J\subset G_1,\ J=<a_1,a_2,x_1,x_2>$.
Then $H/J$ is a two-sided totally geodesic hypersurface in $M_1$ and 
we may cut $M_1$ along this hypersurface to get a connected manifold $M'_1$ 
that has two boundary components which are isometric 3-dimensional hyperbolic
manifolds given by $H/J$.  We may do the same with $M_2$ to get $M'_2$ which is
also connected (the subgroup $J$ in question doesn't change). 

\midinsert
$$\hss\psfig{file=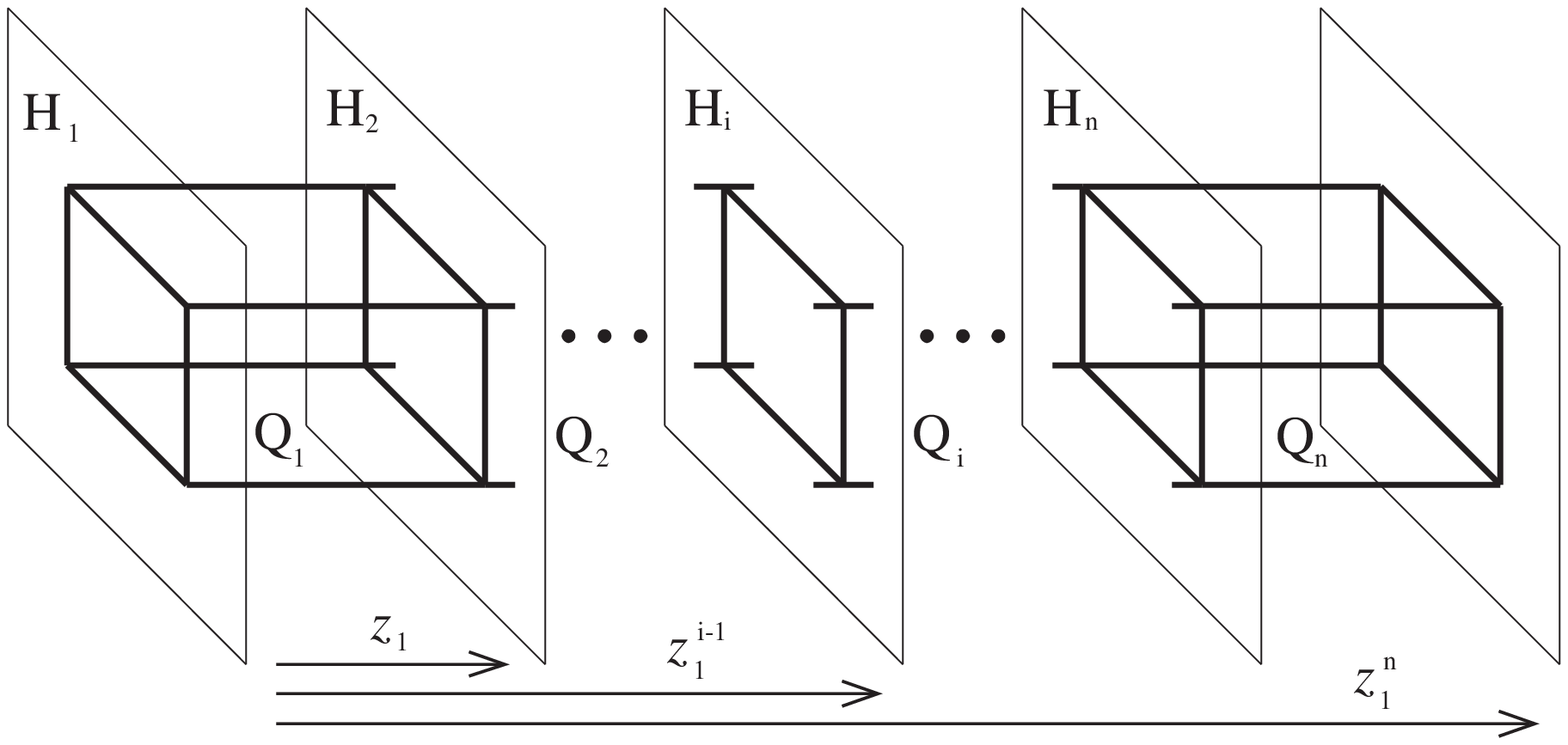,height=2.5truein}\hss$$
\botcaption{Figure 10}
The planes $H_i$ in the polyhedron $Q$
\endcaption
\endinsert

Now let $M=\hiv/G$, where $G$ is a subgroup generated by
any of the side-pairings of the polyhedron $Q$
defined in the previous section.  We identify the block $Q_1$ with $P$
and its side-pairing with $\fai$ or $\faii$.
Let $H_i$ be the hyperplane supporting the side $Z_1$ of the block $Q_i$,
$i=1,\dots,n$.  So, $H_i=z_1^{i-1}(H)$ and it is kept invariant by the
subgroup $J_i\subset G$, $J_i=z_1^{i-1}J z_1^{-(i-1)}$ (see Figure 10).
Suppose we know that $H_i$ is precisely invariant under the subgroup
$J_i\subset G$ . Then we may cut $M$ along the hypersurfaces
$H_i/J_i$.
The hyperplanes $H_i$ are exactly the ones that separate the
polyhedron $Q$ into blocks $Q_1,\dots,Q_n$ (see Figure 10).  Under
identification of
paired sides, we see that each block yields a submanifold of type $M'_1$
or $M'_2$.  Therefore, after cutting along $H_i/J_i,\ i=1,\dots,n$,
we will get $n$ pieces, each isometric to either
$M'_1$ or $M'_2$.  Thus, we will have shown

\proclaim{Proposition 4.1}
Any of the manifolds constructed in section 3 are obtained by gluing $n$
copies of either $M_1'$ or $M_2'$ so that
each copy of $M_1'$ or $M_2'$ is glued to another copy
along one of their two totally geodesic boundaries
i.e. they are strung together in a circular fashion.
Here $M_1'$ and $M_2'$ are manifolds that we obtain by
cutting hyperbolic 4-manifolds $M_1$ and $M_2$ along a totally
geodesic hypersurface.
\qed\endproclaim

\remark{Remark 4.2}
Let $G_{0i},\ i=1,2$ be the subgroup generated by the same generators as 
$G_i$, but with $z_1$ omitted.  Then, by using Maskit's combination theorems
(see theorems 6.19 and 6.24 in \cite{A2}) and fairly standard arguments
(see, for example, \cite{A3}) it is possible to show that
$$
G=\left(G_{0i_1}\underset{J_2}\to\ast
(z_1G_{0i_2}z_1^{-1})\underset{J_3}\to\ast
(z_1^2 G_{0i_3}z_1^{-2})\underset{J_4}\to\ast\dots
\underset{J_n}\to\ast
(z_1^{n-1}G_{0i_n}z_1^{-(n-1)})\right)\underset{z_1^n}\to\ast.
$$
Here $\underset{J_i}\to\ast$ denotes the free product with amalgamation,
while the last $\underset{z_1^n}\to\ast$ is the HNN-extension of the free
amalgamated product in the parentheses by $z_1^n$, where ${z_1^n}$
conjugates subgroups $J_1$ and $z_1^n J_1 z_1^{-n}$.  The index $i_k$
is 1 or 2 depending on whether the block $Q_k$ is of type $\fai$ or
$\faii$.
\endremark

\remark{Remark 4.3}  
The hyperbolic 3-manifold
$H/J$ has been described by Wielenberg, see Example 3 in \cite{Wi2}.
It is the complement of a certain four-component link in $S^3$.  Thus, we
are gluing 4-manifolds along boundaries that are link complements in $S^3$.
\endremark
\vskip6pt

The only thing left to do is to verify that the hyperplanes 
$H,H_1,\dots,H_n$ are precisely invariant, respectively,
under the subgroups $J,J_1,\dots,J_n\subset G$.  To this end we
use the technical theorem stated below.

\proclaim{Theorem 4.4}
Let $P$ be a fundamental polyhedron for a discrete group
$G\subset\isom(\hn)$ that is generated by some side-pairing of $P$.
Let $H$ be a hyperplane in $\hn$ so that 
$\intr_H (H\cap P)\ne\emptyset$ 
and let $J$ be a subgroup of $G$ that keeps $H$ invariant.
Assume the following three conditions hold: 
\roster
\item $H\cap P$ is a fundamental polyhedron for the action of
$J$ in $H$.
\item If $H$ contains a side $S$ of $P$ and $s$ is the side-pairing 
corresponding to $S$ then $s(H)\ne H$.  
\item Suppose $H$ contains an edge $E$ of $P$.  Let 
$\{\sigma_i=(E_i,S_i,R_i,g_i)\}_{i=1,2\dots}$ be the sequence obtained
by the edge-chase corresponding to $E$ as in the edge cycle condition.
(That is, $E_1=E$, $S_i$ and $R_i$ are the sides that determine $E_i$,
$E_i=S_i\cap R_i$, $g_i(E_i)=E_{i+1}$ and $g_i$ is the
side-pairing that pairs $R_i$ and $S_{i+1}$.)
Let $\alpha$ be the angle between $H$ and $S_1$ and let $\theta_i$ be the
dihedral angle of $P$ at the edge $E_i$.  If another edge $E_{l+1}$
in the cycle of $E$ is contained in $H$ and $\beta$ is the angle between
$H$ and $S_l$, then $\theta_1+\dots+\theta_l+\beta-\alpha =k\pi$ must
be satisfied for some integer $k$.
\endroster
Then $H$ is precisely invariant under $J$.
\endproclaim

\demo{Proof}
Let $f\in G$ and suppose that $K=f(H)\cap H\ne\emptyset$.   
We want to show that $f\in J$. There are three cases depending on how $K$
intersects the elements of the tiling $\{g(P),\ g\in G\}$.
We will repeatedly use the fact that 
$H\subset\cup_{j\in J}\, j(P)$ and $f(H)\subset\cup_{j\in J}\, fj(P)$
which follows from the assumption that
$H\cap P$ is a fundamental polyhedron for $J$ in $H$.

{\it Case 1.} There exists an $x\in K$ so that $x\in\intr g(P)$ for some
$g\in G$.
Since $x\in H$, there must be a translate of $P$ under $J$ that contains
$x$.  The only possible candidate is $g$, so we conclude $g\in J$.
Likewise, since $x\in f(H)$, there is a $j\in J$ so that $x\in fj(P)$.
But this can only happen if $fj=g$, so $f=gj^{-1}\in J$.

{\it Case 2.} There exists an $x\in K$ that is contained in the interior of
a side $R$ of some translate of $P$.  Then $R$ is common to exactly two
translates of $P$. If $R$ is not contained in $H$ then $H$ cuts into the
interior of both of those translates.  Again, the parts of $H$ that are in
the interiors of these translates must be covered by translates of $P$
under $J$ so the translates abutting $R$ are of form $j(P)$ and $j'(P)$
for some $j,j'\in J$.  Furthermore, there must be a
$j''\in J$ so that $x\in fj''(P)$.  Since $x$ is in only two translates of
$P$, this means that either $fj''=j$ or $fj''=j'$.  In both cases we get
$f\in J$.  

If, on the other hand, $H$ does contain $R$, then at least one
of the two translates of $P$ abutting $R$ is of form $j(P),\ j\in J$.  
Assuming  $f(H)\ne H$ gives us that $f(H)$ intersects the interior of
$j(P)$.  The portion of $f(H)$ in $\intr j(P)$ must be in some $fj'(P)$
for some $j'\in J$ so we get $fj'=j$, which forces $f\in J$, contradicting
$f(H)\ne H$.  Therefore $f(H)=H$.  Now the two translates of
$P$ that abut $R$ are of form $j(P)$ and $js^{-1}(P)$, where $s$ is the
side-pairing of the side $S$ for which $j(S)=R$.  One of those translates is
also of form $fj'(P)$ for some $j'\in J$.  If $fj'=j$ then $f\in
J$.  The other case, $fj'=js^{-1}$ implies $s=j'{} ^{-1}f^{-1}j$, so
$s$ preserves $H$.  This, however, contradicts assumption (2), because
$H$ contains $S$, since $R\subset H$, $S=j^{-1}(R)$ and $j^{-1}(H)=H$.

{\it Case 3.} If neither case 1 nor 2 occurs, we get that $K$ is contained
in translates of edges of $P$, which are $n-2$-dimensional.  Since 
$\dim K\ge n-2$ we get that $K$ must be $(n-2)$-dimensional, which implies
$f(H)\ne H$.
Furthermore, there exists an $x\in K$ and an edge $E'$ of some $g(P)$ so
that $x$ is in the interior of $E'$.  As before, one of the translates of
$P$ that contains $x$ must be of the form $j(P)$.  Move everything by
$j^{-1}$ so that $x$
is now on an edge $E$ of $P$ and $E\subset j^{-1}f(H)\cap H$.  The
translates of $P$ that abut $E$ are 
$P,\ g_1^{-1}(P),\ g_1^{-1}g_2^{-1}(P),\dots$, so as before, there must be
a $j'\in J$ and an integer $l$ so that
$j^{-1}fj'=g_1^{-1}\circ\dots\circ g_l^{-1}$.
But then $j^{-1}f(H)=j^{-1}fj'(H)=g_1^{-1}\circ\dots\circ g_l^{-1}(H)$, so 
$E\subset g_1^{-1}\circ\dots\circ g_l^{-1}(H)\cap H$.
From $E\subset g_1^{-1}\circ\dots\circ g_l^{-1}(H)$ we get that 
$E_{l+1}=g_l\circ\dots\circ g_1(E) \subset H$, so $E_{l+1}$ is in the
cycle of $E$ and is contained in $H$.
Let $K^\perp$ be the 2-dimensional orthogonal complement of $K$ through $x$.  
The intersections of translates of $P$ that abut $E$ with $K^\perp$
are angles with rays emanating from a single vertex $x$.  Intersections of
$H$ and  $g_1^{-1}\circ\dots\circ g_l^{-1}(H)$ with $K^\perp$ are two
lines and the angle between them is
$\theta_1+\dots+\theta_l+\beta-\alpha$.
Condition 2) now says that this angle is $k\pi$, so the lines are
identical and so are the hyperplanes that they represent.
From here it follows that $f(H)=H$, 
a contradiction with $f(H)\ne H$.  Therefore, case 3 never occurs and
$f\in J$ by cases 1 and 2.
\qed\enddemo

\remark{Remark 4.5} The proof of Theorem 4.4 did not use any hyperbolic
space-specific properties, only the fact that $P$ was a fundamental
polyhedron.  Therefore it also applies in the other two constant
curvature settings, that is, for fundamental polyhedra of discrete isometry
groups of the $n$-sphere and Euclidean $n$-space.
\endremark

\remark{Remark 4.6} Notice that the group $G$ in the theorem did not have to
be torsion-free.  However, if $H$ contains an edge of $P$,
condition (3) allows the number $k$ that was defined in the edge cycle
condition to only be 1 or 2.
\endremark
\vskip6pt

\example{Examples} We give several applications of the theorem that include
the claims of precise invariantness needed for Proposition 4.2.
All except example 4.8 have as $P$ the polyhedron defined in section 1.
\endexample

\remark{Example 4.7}
Let $G=G_1$ or $G_2$, $H$=the hyperplane based on the plane
$\{z=0\}$, $J=<x_1,y_1,a_3,a_4>$.  Clearly $H$ is invariant under $J$.
By applying (now in dimension 3) Poincar\' e's polyhedron theorem to
$H\cap P$ and restrictions of $x_1,y_1,a_3,a_4$ to $H$ we may easily see that
$H\cap P$ is a fundamental polyhedron for $J$ in $H$.  (Here conditions
(2) and (3) from the theorem do not apply.)  Therefore, $H/J$ is a
totally geodesic hypersurface embedded in $M_1$ or $M_2$.
\endremark
\vskip6pt

\remark{Example 4.8}
Let $G=G_1$ or $G_2$, $H=Z_1$ and $J=<x_1,y_1,a_1,a_2>$.  As in Example 4.7
we check that $H\cap P$ is a fundamental polyhedron for $J$ in $H$.
Here we also need to verify condition (3) of Theorem 4.4. (Condition (2)
clearly holds.) Taking, for
example, $E=Z_1 \cap A_1$ whose cycle is
$\{Z_1\cap A_1,\ A_1'\cap Z_1,\ Z_1'\cap A_5,\ A_5'\cap Z_1', \}$ we see that
$\alpha=0,\ l=1$ and $\beta=\pi/2$, so condition (3) is satisfied.
Using Theorem 4.4 gives us that $H/J$ is a
totally geodesic hypersurface embedded in $M_1$ or $M_2$.
\endremark
\vskip6pt

\remark{Example 4.9} It is now easy to see that the hyperplanes 
$H_1,\dots,H_n$ (in above
notation) are precisely invariant under the subgroups
$J_1,\dots,J_n\subset G$.
The proof for $H_1$ and $J_1$ corresponds to the one in Example 4.8,
while the other cases correspond to Example 4.7.  This completes the proof of
Proposition 4.2.
\endremark
\vskip6pt

\remark{Example 4.10}
Let $G=G_1$ or $G_2$ and let $H_1,\ H_2$ be the hyperplanes
based respectively
on the planes $\{x-y=0\}$ and $\{x+y=0\}$.  We may use Theorem 4.4 to
verify that $H_1$ is precisely invariant under
$<a_2,a_4,a_6,b_1,d_2,d_4,y_1x_1>$ and that $H_2$ is precisely invariant
under $<a_1,a_3,a_5,b_1,d_1,d_3,y_1^{-1}x_1>$.  Again, condition (1)
of Theorem 4.4 is verified by the Poincar\'e polyhedron theorem in
dimension 3.  Note that condition (3) of that same theorem applies.
\endremark

\head
5. Proof of Theorem B
\endhead

We will show that the polyhedron $P$ has volume $2\cdot \const$.  Then
every manifold obtained from any of the side-pairings of $Q$ described
in section 3 will have volume $2n\cdot\const$.
The Gauss-Bonnet formula (see \cite{G}, page 84, and \cite{H}) applied to a
non-compact hyperbolic
4-manifold $M$ gives $\volm M=\chi(M)\cdot\const$. Here $\chi(M)$ is the
Euler characteristic of the compact part of $M$, i.e. the manifold with
boundary obtained by retracting every end $E\times [0,\infty)$ of $M$ to
$E\times \{0\}$.  It will be enough to show that
the manifold obtained from either of the side-pairings $\fai$ or $\faii$ has
Euler characteristic 2.

Recall that for a finite CW-complex $X$, $\chi(X)$ may be computed
either as an alternating sum of the numbers of $i$-cells in $X$, or as
the alternating sum of the ranks of the $i$-th homology groups of $X$.

Let $P$ be finite-sided $n$-dimensional hyperbolic polyhedron with a
side-pairing defined on it.  Now
let $X$ be the CW-complex obtained from $P$ in the obvious way, with
0-cells the real and ideal vertices of $P$, 1-cells the 1-sides of
$P$ together with their points at infinity, and so on.  Then $X$ inherits
identifications by side-pairings of
$P$, which give rise to a quotient space $Y$, also a CW-complex (even
if $P$ yields a manifold by identification, $Y$ will not be one).  We then
have

\proclaim{Lemma 5.1}
If $M$ is obtained from a side-pairing of $P$ and $Y$ is as above, then
$$\chi(M)=\chi(Y)-\text{number of ends of $M$}.$$
In other words, we may compute $\chi(M)$ directly from the 
polyhedron by taking the alternating sum of numbers of cycles of 
$i$-sides and ignoring the cycles of ideal vertices.
\endproclaim

\demo{Proof}
Using $M$ to also denote the compact part of the hyperbolic manifold, we
see that $Y=M\cup V$, where $V$ is a disjoint union of cones over
Euclidean manifolds that are the boundary of $M$, so that $M\cap V$ is a
disjoint union of Euclidean manifolds.  Consider the absolute Mayer-Vietoris
sequence for $M$ and $V$ (see \cite{Do}, Proposition 8.15).  Enumerate the
terms so that the $k$-th
homologies of $M\cup V$, the sum of $M$ and $V$, and $M\cap V$ correspond
to indices $3k,\ 3k+1$ and $3k+2$ respectively ($k\ge0$).
Let $c_j$ and $z_j$ denote respectively
the rank of the $j$-th term and the rank of the kernel of the 
homomorphism joining the $j$-th and the $(j-1)$-st terms of the sequence.
By exactness of that sequence we have
$c_j=z_{j-1}+z_j$.  Use this equality to see that
$\sum (-1)^k c_{3k}-\sum (-1)^k c_{3k+1}+\sum (-1)^k c_{3k+2}=0$, which
is exactly  \linebreak
$\chi(M\cup V)-(\chi(M)+\chi(V))+\chi(M\cap V)=0$.  Since $V$
is contractible, $\chi(V)=$number of components of $V$=number of
of boundary components of $M$.  The fact that $M\cap V$ is a disjoint
union of Euclidean manifolds
implies $\chi(M\cap V)=0$ which yields the desired formula.
\qed\enddemo

Now to finish the proof of Theorem B, we just have to count cycles of
$i$-sides for $\fai$ or $\faii$.  The polyhedron
$P$ from section 1 has one 4-side, 36 3-sides, 168 2-sides, 216 1-sides and
48 real 0-sides.  Each 3-side is paired to exactly one other one, which yields 18
cycles of 3-sides.  Among 2-sides, there are 24 with dihedral angle
$\pi/4$, giving 3 cycles, and 144 with dihedral angle $\pi/2$,
yielding 36 cycles.  Among 1-sides, there are 80 with normalized solid angle
1/16 giving 5 cycles and 136 with normalized solid angle 1/8, yielding
17 cycles.  From Figure 6 we know there are 2 cycles of 0-vertices.
Thus, $\chi(M)=2-(17+5)+(36+3)-18+1=2$, which completes the proof.
\qed

\remark{Remark 5.2}
Note that the same reasoning as in the above paragraph may be used to
see that, when the side-pairing of an $n$-polyhedron $P$ yields a
manifold $M$, $\chi(M)$ depends only on the alternating sum of normalized
solid angles of $P$ over all the $i$-sides of $P$.  In particular, it
doesn't depend on the side-pairing of $P$.  Indeed, the sum of
normalized solid angles for each cycle of an $i$-side is exactly 1,
which is contributed to the count of cycles of $i$-sides.
\endremark

\enddocument

\end